\baselineskip=18pt

\font\bigbf=cmbx10 scaled\magstep2

\def\proclaim #1. #2\par{\medbreak
\noindent{\bf#1. \enspace}{\sl#2}\par\medbreak}

\def\section #1. #2\par{\bigbreak
\noindent{\bf#1. \enspace}{\bf#2}\par\medbreak} 

\def\qed{{\vrule height 5pt width 5pt depth 0pt}} 
\def\t{\varphi}
\def\q0{q_0}
\def\p0{\t_0}
\def\S{{\bf S}}

\def\R{{\bf R}}
\def\e{\varepsilon}
\def\Z{{\bf Z}}
\def\M{{\cal M}}

\def\qp{q^{\prime}}

\def\Rt{R_{\theta}}
\def\qp{q^{\prime}}

\def\At{{\tilde A}}
\centerline{\bigbf Some Remarks on Periodic Billiard Orbits}
\centerline{\bigbf in Rational Polygons}
\medskip
\centerline{M. Boshernitzan\footnote{$^{(1)}$}{Department of Mathematics, Rice
University}, G. Galperin\footnote{$^{(2)}$}{Forschungszentrum BiBos, 
Universit\"at Bielefeld}\footnote{$^{(3)}$}{Current address: Department
of Mathematics, Eastern Illinois University}, 
T. Kr\"uger$^{(2)},$ and 
S. Troubetzkoy$^{(2)}$\footnote{$^{(4)}$}{Institute for Mathematical Science, 
SUNY Stony Brook}}

\section 1. Introduction
\par
A billiard ball, i.e.~a point mass,
moves inside a polygon $Q$ with unit speed 
along a straight line until it reaches the boundary $\partial Q$ of the 
polygon, then 
instantaneously changes direction according to the mirror law: ``the angle of 
incidence is equal to the angle of reflection,'' and continues along the 
new line (Fig.~1(a)).  
Despite the simplicity of this description there is much that is unknown
about the existence and the description of periodic orbits in arbitrary
polygons. On the other hand, quite a bit is known about a special class of
polygons, namely, rational polygons. A polygon is called
{\it rational} if the angle between each pair of sides is a rational 
multiple of $\pi.$  
The main theorem we will prove is

\proclaim Theorem 1.  For rational polygons, periodic points of the billiard 
flow are dense in the 
phase space $\M$ of the billiard flow.\par

Theorem 1 is a strengthening of Masur's theorem [M] who has shown that any
rational polygon has ``many'' periodic billiard trajectories; more
precisely, the set of directions of the periodic trajectories are dense
in the set of velocity directions $\S^1.$ We will also prove some
refinements of theorem 1: the ``well distribution''
of periodic orbits in the polygon and the residuality of the points
$q \in Q$ with a dense set of periodic directions (precise statements
of these results will be given in section 3).

The structure of the article is as follows.  In section 2 we will
give a brief description of billiards in polygons and some results
related to theorem 1.  In section 3 we will state the strengthenings
of theorem 1, and the proofs will be given in section 4.

\section 2. Description of billiards in polygons
\par
The trace of a moving billiard ball is called a {\it billiard trajectory} 
or {\it orbit}.
If a billiard trajectory hits a vertex of the polygon, then it is 
called {\it singular}. For convenience we will define such billiard 
trajectories by 
continuity from the left (with respect to a fixed orientation of the boundary),
thus every trajectory is infinite and
the singular trajectories are the discontinuities of the flow.  
It is convenient to consider the set of singular trajectory segments 
which start at a vertex and end at a vertex.  Such segments are called
{\it generalized diagonals} and the number of links is called the {\it 
length}
of the generalized diagonal (Fig.~1(b)).  We remark that the set 
of generalized diagonals is countable. 
If a billiard ball returns to its initial 
position and has the same velocity direction, then its orbit is
called {\it periodic} (Fig.~1(c)).

There is a very useful tool in the analysis
of billiards in polygons: the {\it unfolding} of trajectories.
Instead of reflecting the trajectory with respect to the side of the 
polygon it hits we reflect the polygon itself with respect to the same
side.  The two adjacent
links of the trajectory then become part of a straight line
(Fig.~2(a)).  Continuing this procedure for ever unfolds the trajectory into
a half-line through a forward {\it corridor} of polygons.
The backward trajectory can be similarly unfolded.
We enumerate the unfolded polygons starting with zero, while the
length of an unfolding is the number of polygons it contains (Fig.~2(b)). 

The set of pairs $\big \{ (q,v) \big \},$ where $q \in Q$ is a position of
the ball and $v \in \S^1$ is its velocity constitute the phase space
$\M$ of the billiard flow.  The phase space $\M$ can be thought of as a 
right prism foliated with a
collection of ``floors'' $Q(\t) := Q \times \t$ with $\t \in \S^1.$
Then the billiard flow in the phase space can be imagined as a straight line
flow on the floors of the phase space.  When the flow reaches the boundary 
it jumps from the point $(q,\t)$ to the point $(q,\t^{\prime})$ on the floor 
$Q(\t^{\prime})$, where $\t$ 
and $\t^{\prime}$ are related by the mirror law (Fig.~3).  In the case 
the polygon is rational, the flow is restricted to a finite number of 
floors.  Using the identifications of the boundary one 
can glue these floors
together to produce an orientable surface $R_{\t}$  
(see the survey article
[Gu] for a full explanation).

For any $\t \in S^1$ for which $R_{\t}$ contains no generalized diagonal
the flow $\phi_t$ restricted to the
invariant surface $R_{\t}$ is {\it minimal}, i.e.~the orbit of every point is
dense [Gu]. 
Using Teichm\"uller theory Kerckhoff, Masur, and Smillie showed that for
almost every $\t \in \S^1$ the flow $\phi_t|_{R_{\t}}$ is {\it uniquely 
ergodic}, i.e.~the only ergodic invariant measure is the
Lebesgue measure [KMS]. Using similar techniques
Masur has shown that for a dense set of directions $\t \in S^1$ 
the flow $\phi_t$ restricted to $R_{\t}$ has at least one periodic point [M].  
However, his result gives no indication of the location of the periodic orbit
on the invariant surface $R_{\t}.$

\section 3. Results
\par

Throughout this section we will assume that $Q$ is a rational polygon.

\proclaim Theorem 1. Periodic points of the billiard flow are 
dense in the phase space $\M$.\par

We will also prove a slightly stronger theorem.  
For this purpose we need one more definition. 
A periodic orbit $\gamma$ is called $\e$-{\it well distributed} on the table if
for every convex set $A \subset Q$ 
$$\left | {{\hbox{length}(\gamma \cap A)} \over {\hbox{length}(\gamma)}} - 
\mu_{\scriptstyle {\rm Leb}}(A) \right
| < \e.\eqno(1)$$
Here the length of a periodic orbit is its geometric length, not the 
number of links in its trajectory.

\proclaim Theorem 2.  The set of $\e$-well distributed periodic points
of the billiard flow  is dense in the phase space $\M$ for every $\e > 0$.\par

Let  $V$ be the set of $Q$'s vertices and
$$G := \Big \{q \in Q \backslash V: \, \forall \e > 0, 
\hbox{ for a dense set of directions } \t, \,
(q,\t) \hbox{ is an } 
\e\hbox{-well distributed periodic point} \Big \}.$$
As another refinement of theorem 1 we have: 

\medbreak\noindent
{\bf Theorem 3.}\enspace{\sl $G$ is residual.}\footnote{$^{(1)}$}{A set is 
called 
{\it residual} if it contains a dense $G_{\delta}$ set, i.e.~a countable 
intersection of open dense sets.}
\par\medbreak

\noindent The natural question arises: does $G = Q \backslash V$ for all
rational polygons? Does $G = Q \backslash V$ mod(0) for all polygons?

Equality in the stronger sense holds for Coxeter chambers 
(polygons which tile the plane under their reflection group), almost 
integrable polygons (finite unions of Coxeter chambers) [Gu], 
regular polygons, and closely related polygons [V].  Next we consider
$$B := \Big \{ q \in Q \backslash V: \hbox{ there is no 
direction } \t \hbox{ for which } (q,\t) \hbox{ is periodic} \Big \}.$$

We remind the reader that we have a standing assumption that $Q$ is a
rational polygon. 

\proclaim Theorem 4. If $Q$ is convex then $B$ is contained in a finite union 
of segments.  If 
$Q$ is a triangle, then $B$ is at most finite.\par

Finally we turn to the question: for which angles $\t$
is there at least one orbit which is not $\e$-dense?
We show:

\proclaim Theorem 5. The set of directions for which there exists a non 
$\e$-dense orbit is an at most countable closed set for any
$\e > 0.$ There are examples when this set is not finite.\par

\section 4. Proofs of theorems
\par
\noindent
{\bf Proof of theorem 1:}

First of all we fix a uniquely ergodic direction $\theta \in \S^1.$  
As mentioned above by avoiding a countable set of $\theta$'s we can assume that 
$\Rt$ contains no generalized diagonals.
We claim that using the unique ergodicity we can choose 
$N$ so large that for all $x \in \Rt$
the first $N$ links of $x$'s-forward orbit and the first $N$ 
links of $x$'s-backward orbit are both $\e/2$-dense 
in the surface $\Rt.$\footnote{$^{(2)}$}{This fact can also be derived from the
minimality of the direction $\theta.$}
To see this note that the unique ergodicity of the 
flow implies that for every continuous function $g$ the ergodic average
$1/T \int_0^T g(\phi_tx)\,dt$ converges uniformly to
$\int_{\M}g\, d\mu_{\rm Leb}$ (see [W] for the proof
which holds without change in this more general setting).

There are only a finite number of generalized diagonals of length less
than or equal to $2N.$  In the phase space $\M$ these generalized diagonals
lie on a finite number of floors $Q(\t_1), Q(\t_2), \dots, Q(\t_{k(N)}).$
Let $\delta := \delta_N > 0$ be so small that for 
any $\theta^{\prime}$ satisfying
$|\theta - \theta^{\prime}| < \delta$ we have $Q(\t_i) \cap 
R_{\theta^{\prime}} = \emptyset$ for $i = 1, \dots, k(N)$ (Fig.~4).
Therefore, if a generalized diagonal belongs to  $R_{\theta^{\prime}},$ 
then its length is greater than $2N.$ 

From Masur's theorem we know that there is a periodic point whose direction 
is arbitrarily close to $\theta.$ In particular, we can find a periodic 
point $(\q0,\p0) \in \M$ satisfying
$|\p0 - \theta| < \delta$ and $< \e/(2 \cdot \hbox{const} \cdot
N \cdot \hbox{diam }Q).$ Here, const is a constant depending on
$\theta,$ which will be defined later. 
We consider the point
$(\q0,\theta)$ and its forward unfolding of length $N$ and its backward
unfolding of length $N$.
We claim that {\it either} the corridor of length $N$ for 
the forward trajectory of $(\q0,\p0)$
{\it or} the corridor of length $N$ for 
the backward trajectory of $(\q0,\p0)$
coincides with the same length corridor of $(\q0,\theta).$

Suppose this is not true. The forward corridors
of  $(\q0,\theta)$ and $(\q0,\p0)$ coincide for a while.
Define $j_1 < N$ to be the length of the forward part of the corridors of
$(\q0,\theta)$ and $(\q0,\p0)$ that coincide. The two corridors ``branch''
apart at polygon number $j_1 -1$; let $A$ be the common vertex of branching 
corridors.
Similarly, let $j_2 < N$ be the length of the backward corridors and
$B$ the common vertex of the backward branching corridors.
Then the straight line segment $AB$ inside the corridor of $(\q0,\theta
)$ is a generalized diagonal of length $j_1 + j_2 -1 < 2N.$
The direction of $AB$ lies in the interval of directions $\{t :
\theta \le t \le \p0\}$ 
(we assume $\p0 > \theta$) and
thus is its distance to $\theta$ is less than $\delta$ (Fig.~5). 
This contradicts the choice of $\delta.$

Without loss of generality
we will assume that both the periodic trajectory and the uniquely ergodic 
trajectory stay in the same forward corridor.
The endpoints of these two trajectories lie on the boundary of the $(N-1)$st
polygon.
Since $|\p0 - \theta| < \e/(2 \cdot \hbox{const} \cdot N \cdot
\hbox{diam} Q)$  
it follows that if the const is small enough, then 
the distance between the endpoints of the two unfolded 
trajectories  is less than $\e/2$ (Fig.~6).  
Furthermore, the first $N$ links of the uniquely ergodic trajectory are
$\e/2$-dense in $\Rt.$ These two facts combined show that
the periodic trajectory is $\e$-dense in $R_{\p0}$. This can be 
done for every uniquely ergodic direction $\theta$. The set 
$\{ R_{\theta} : \phi|_{R_{\theta}}
\hbox{ is uniquely ergodic}\}$ is dense in the whole phase space $\M.$ 
Since $\e$ is arbitrary
this completes the proof of theorem 1. 
\hfill\qed

\medskip\noindent
{\bf Proof of theorem 2:}

We given an equivalent definition of $\e$-{\it well distribution}.
We fix an embedding $Q \subset \R^2.$
Let $A_{p,q,r,s}$ be the intersection of $Q$ with the open ball with center
$(p/q,r/s)$ and diameter $1/(\max(q,s)).$ 
Enumerating gives rise to the countable basis $\{\At_i\}.$ 
A periodic orbit $\gamma$ is $\e$-{\it well distributed} on the table if
for each set $\At_i$ with diam$(\At_i) > \e$
$$\left | {{\hbox{length}(\gamma \cap \At_i)} \over {\hbox{length}(\gamma)}} - 
\mu_{\scriptstyle {\rm Leb}}(\At_i) \right
| < \e.\eqno(1)$$
For a fixed sufficiently small $\e > 0$ for a convex set $A$ with sufficiently
small diameter equation 1 automatically holds, and 
we can approximate any convex set $A$ with large diameter by finite unions and
intersections of the $\At_j.$
Using this reasoning one can show 
that for every $\delta > 0$ there is a $\e > 0$
so that any $\e$-well distributed point in the sense above is $\delta$-well
distributed in the sense of section 2 and $\e \rightarrow 0$ 
as $\delta \rightarrow 0$.  

We need to modify several steps in the proof of theorem 1. 
In the proof of theorem 1 we choose $N$ so large that
the first $N$ links of $x$'s forward and backward orbit are both 
$\e/2$-dense. 
The proof given there actually shows the stronger result, both orbit segments 
are $\e/2$-well distributed as well (the definition of well distributed is 
analogous to the one for periodic points introduced in section 3).
We choose the integer $N$ slightly larger so that the first $N$ links are 
$\e/3$-well distributed. 
We also choose $N$ so large that for any $k \ge N$ if the first $k$ links 
of the orbit of any
point $x$ are $2\e/3$-well distributed then the first $k+1$ links are
$\e$-well distributed.

We choose $\delta := \delta_{2N}$ and the periodic point $(q_0,\p0)$
so that it satisfies $|\p0 - \theta| < \delta_{2N}$ and 
$< \e/(3 \cdot \hbox{const} \cdot
2N \cdot \hbox{diam}Q).$ 
Then the periodic trajectory of $(q_0,\p0)$
shadows the $\e/3$-well distributed uniquely ergodic trajectory for $2N$
links either forward or backward.  Thus the first $2N$ links of this orbit
are $2\e/3$-well distributed and the first $2N+1$ links are $\e$-well
distributed.  Suppose the number of links of the 
$(q_0,\p0)$ trajectory is $M.$ If $M \le 2N + 1,$ then the trajectory is
clearly $\e$-well distributed. 

If $M > 2N +1$ let the cyclically ordered set
$\{L_1,L_2,\dots,L_M\}$ denote the links of the trajectory.  
We want to {\it partition} the trajectory into $\e$-well distributed trajectory
segments of lengths between $N$ and $2N.$ 
From formula (1) it is clear that the concatenation of 
$\e$-well distributed trajectory segments is itself $\e$-well
distributed. Therefore, 
the construction of such a partition will finish the proof.

To construct such a partition we start by {\it covering} the trajectory by 
trajectory segments, i.e.~by ordered sets $\{ L_i, \dots, L_{i+2N-1}\}$
of length $2N$ that 
are $\e$-well distributed and for which all its ordered subsets $\{L_j,\dots,
L_{j+k-1}\} \subset \{ L_i, \dots, L_{i+2N-1}\}$
of length $k \ge N$ are also $\e$-well distributed.
Without loss of generality we suppose that $\{L_1,\dots
L_{2N}\}$ as well as all its consecutively ordered subsets of length
$k \ge N$ are $2\e/3$-well distributed.  Now we apply this argument again to the
point $\phi_{T(2N+1)}(q_0,\p0)$
to conclude that either $\{L_2,\dots,L_{2N+1}\}$ or
$\{L_{2N+2},\dots,L_{4N+1}\}$ has this property. If the latter occurs, 
then the link $L_{2N+1}$ is not covered. In this case we replace $\{L_1,\dots,
L_{2N}\}$ by $\{L_1,\dots,L_{2N+1}\}$ which is $\e$-well distributed.
Continuing this process
inductively yields the desired covering (Fig.~7).
To finish the proof we split the covering into a partition of  
$\e$-well distributed pieces
of different lengths, but with all the lengths being at least $N$ and at
most $2N.$ \hfill\qed

\medskip\noindent {\bf Proof of theorem 3:}

Fix a uniquely ergodic direction $\theta$.  Let 
$$A^{\delta}_{\e}(\theta) := \Big \{ q \in Q \backslash V : 
(q,\t) \hbox{ is an } \e\hbox{-well distributed periodic point for 
some } \t, \ |\t - \theta| < \delta \Big \}.$$
In the proof of theorem 1 we showed that 
for all uniquely ergodic directions $\theta$ and for all $\e> 0$ 
there is a $\delta := \delta(\theta,\e) > 0$ such that 
$A^{\delta}_{\e}(\theta)$ is dense in $Q.$  The set
$A^{\delta}_{\e}(\theta)$ is also open because each periodic point is
contained in an open strip, that is if $x = (q,\t)$ is periodic, then there 
is an open disc $D \subset Q \backslash V$ such that $(\qp,\t)$ is periodic 
for all $\qp \in D$ (see [GKT]).
Choose a countable dense set 
of uniquely ergodic directions $\{ \theta_i \} \subset \S^1$ and 
$\e_i > 0$ satisfying $\lim_{i \rightarrow \infty} \e_i = 0.$ 
Furthermore choose positive numbers
$\delta_i \le  \delta(\theta_i,\e_i)$ such that 
$\lim_{i \rightarrow \infty} \delta_i = 0.$ Then
$$\bigcap_{i \in \Z^+} A^{\delta_i}_{\e_i}(\theta_i) \subset G,$$
and thus $G$ is residual.\hfill\qed

\bigskip
\noindent {\bf Proof of theorem 4:}

First we will show that the set $B$ 
is contained in a finite number of line segments.
To do this we will consider only billiard trajectories which hit some side 
of the polygon perpendicularly. Let $\t_i$ be the direction perpendicular to 
the $i$th side of $Q.$ As discussed in section 2, the invariant surface
$R_{\t_i}$ contains at most a finite number of generalized diagonals.
Any perpendicular 
trajectory which enters a vertex is automatically a generalized diagonal 
because its backward orbit also enters the same vertex. This means that
only a finite number of points on the $i$th side have singular perpendicular 
trajectories. The number of such diagonals is clearly
less than the number of floors of the invariant surface $R_{\t_i}$
multiplied by the number of vertices of $Q$. We remark that the number of 
$R_{\t_i}$'s floors 
is less than the greatest common denominator of all the inner angles of 
the polygon.

The perpendicular orbit through any other 
point on the $i$th side hits that side perpendicularly twice and thus
is periodic (for details see [Bo][GSV]).
For each point $q \in Q$ we consider the perpendiculars to the
straight lines containing the sides of $Q$. Since $Q$ is convex, for 
the nearest side of $Q,$
the base point of this perpendicular will belong to an interior point of
the side and not to its continuation. Consequently, any point $q$ which
is not covered by one of the finite number of generalized diagonals
from the surfaces $\cup_i R_{\t_i}$ has a ``perpendicular'' periodic orbit
passing through it.

We remark that using the ideas from the proofs of theorems 1-3 we can show
a slightly stronger result. Namely, $B$ is contained in a Cantor subset 
of these segments.

If $Q$ is an acute or right triangle, then for each
$q \in Q \backslash V$ the perpendicular to the lines containing the sides lie 
inside $Q$ (if $q \in \partial Q,$ then this perpendicular is just the 
point $q$).
For obtuse triangles the same is true for at least two sides.
Thus for rational triangles, for each $q \in Q \backslash V$ 
there are at least two distinct singular perpendicular billiard
orbits passing through $q.$
Each point $q \in Q \backslash V$ which is not covered by a perpendicular
periodic trajectory must be covered by two 
distinct perpendicular generalized diagonals. However, any two generalized 
diagonals which intersect do so transversely or else they would coincide.
If $d_1,\dots,d_k$ denote the perpendicular generalized diagonals, then 
$B$ is contained in $\cup_{i \ne j} (d_i \cap d_j ),$ a finite set.
\hfill\qed

{\it Question:} For which polygons are all $q \in Q \backslash V$ at least
double covered by perpendicular trajectories?  
Note that this does not hold for the regular hexagon although in this case it
is easy to see that $B$ is empty.
\bigskip
\noindent {\bf Proof of theorem 5:}

Let 
$$C_{\e} :=
\Big \{ \t: \exists x \in R_{\t} \hbox{ such that the orbit of }
x \hbox{ is not } \e\hbox{-dense } \Big \}.$$
If $R_{\t}$ contains no generalized diagonal, then the flow
$\phi|_{R_{\t}}$ is minimal
and thus every billiard orbit is dense.  The set of generalized diagonals is at
most countable and thus so is $C_{\e}.$  

The set $C_{\e}$ is clearly closed, for if $\t_i \in C_{\e}$ and $x_i \in
R_{\t_i}$ is not an $\e$-well distributed point, then any weak limit point
of the $x_i$ is also not an $\e$-well distributed point.

An example where $C_{\e}$ is countable is the following. We consider the
L-shaped figure consisting of three identical squares 
as depicted in figure 8(a).
In this figure a periodic orbit with 4 links which avoids the right hand
square is shown.  In figure 8(b), an unfolding of length 3 of the L-shaped
figure along with a periodic orbit with 8 links is drawn. This periodic
orbit also never enters the right hand square.  In 
general, the analogous unfoldings of length
$k$ contains a periodic orbit with $2k+2$ links which avoids 
the right hand square.
\hfill\qed

{\it Question:} are there convex examples of 
this phenomenon?

\section 5. Acknowledgements
\par
MB is partially supported by NSF-DMS-9224667. 
GG is supported by the Alexander von Humboldt Stiftung.
ST was supported during part of this research by the Deutsche 
Forschungsgemeinschaft. 

\section 6. References
\par
\frenchspacing

\item{[BKM]} C.~Boldrighini, M.~Keane, and F.~Marchetti {\it Billiards in 
polygons} Ann. Prob. 6 (1978) 532-540.

\item{[Bo]} M.~Boshernitzan {\it Billiards and rational periodic directions in
polygons} Amer.~Math.~Monthly June-July (1992) 522-529.

\item{[GKT]} G.~Galperin, T.~Kr\"uger, and S.~Troubetzkoy
{\it Local instability of orbits in polygonal and polyhedral billiards,} 
Comm.~Math.~Phys.~(1994) (to appear).

\item{[GSV]} G.A.~Galperin, A.M.~Stepin, and Ya.B.~Vorobetz {\it Periodic 
billiard trajectories in polygons: generating mechanisms}
Russian Math.~Surv.~47:3 (1992) 5-80.

\item{[Gu]}  E.~Gutkin {\it Billiards in polygons}  Physica D 19, 
(1986) 311-333.

\item{[KMS]} S.~Kerckhoff, H.~Masur, and J.~Smillie {\it Ergodicity of billiard
flows and quadratic differentials} Annals Math.~124 (1986) 293-311.

\item{[M]} H.~Masur {\it Closed trajectories for quadratic differentials
with an application to billiards} Duke Math.~J.~53 (1986) 307-313.

\item{[V]} W.~Veech {\it The billiard in a regular polygon} Geom.~Func.~Anal.~2
(1992) 341-379.

\item{[W]} P.~Walters {\it An introduction to ergodic theory}, Springer-Verlag,
1982,

\section 7. Figure captions\par

\item{Fig.~1} (a) The mirror law, (b) a generalized diagonal, (c) a 
periodic orbit

\item{Fig.~2} Unfolding a trajectory 

\item{Fig.~3} The mirror law in phase space  

\item{Fig.~4} The $\delta$-neighborhood of $R_{\theta}$ does not contain a 
generalized diagonal of length less than $2N$

\item{Fig.~5} Branching corridors

\item{Fig.~6} $\varepsilon/2$-shadowing of length $N$

\item{Fig.~7} Covering the trajectory by well-distributed trajectory segments

\item{Fig.~8} Periodic orbits in the L-shaped figure which do not enter
the right hand square
\end

\def\etib{\{\bar \eta_t \}_0 ^{\lower 2 pt \hbox{${\scriptstyle \infty}$}}}
\footnote{$^{(2)}$}{i.e.~the only ergodic invariant measure is the
Lebesgue measure. In particular, this means that there is a uniform estimate 
on the time need for all points in the invariant surface to become 
$\e$-well distributed or $\e$-dense.}